\documentclass[numsec,webpdf,modern,medium,namedate]{oup-authoring-template}
\usepackage[T1]{fontenc}
\usepackage[nopatch=footnote, protrusion=false]{microtype}
\usepackage{graphicx}
\usepackage{placeins}
\usepackage{natbib}
\PassOptionsToPackage{numbers, compress}{natbib}

\onecolumn

\graphicspath{{Fig/}}

\usepackage{float} 

\newtheorem{example}{Example}%
\newtheorem{remark}{Remark}%
\usepackage[doublespacing]{setspace}
\usepackage[fontsize=12pt]{fontsize}
\makeatletter
\RequirePackage{xkeyval}
\RequirePackage{tikz}
\RequirePackage{amssymb}

\define@key{boxedtheorem}{titlecolor}{\def\titlecolor{#1}}
\define@key{boxedtheorem}{titlebackground}{\def\titlebackground{#1}}
\define@key{boxedtheorem}{background}{\def\background{#1}}
\define@key{boxedtheorem}{titleboxcolor}{\def\titleboxcolor{#1}}
\define@key{boxedtheorem}{boxcolor}{\def\boxcolor{#1}}
\define@key{boxedtheorem}{thcounter}{\def\thcounter{#1}}
\define@key{boxedtheorem}{size}{\def\size{#1}}
\presetkeys{boxedtheorem}{titlecolor = black, titlebackground = white, background = white,%
                         titleboxcolor = black, boxcolor = black, thcounter=, size = .9\textwidth}{}

\newcommand{\couleurs}[1][]{%
    \setkeys{boxedtheorem}{#1}
    \tikzstyle{fancytitle} =[draw=\titleboxcolor, rounded corners, fill=\titlebackground,
                            text= \titlecolor]
    \tikzstyle{mybox} = [draw=\boxcolor, fill=\background, very thick,
                        rectangle, rounded corners, inner sep=10pt, inner ysep=20pt]
}

\newsavebox{\boiboite}
\newcommand{\titre}{Titre}
\newenvironment{boite}[2][]%
    {%
    \renewcommand{\titre}{#2}
    \couleurs[#1]
    \begin{lrbox}{\boiboite}%
     \begin{minipage}[!h]{\size}
    }%
    {%
     \end{minipage}
    \end{lrbox}
    \begin{center}
    \begin{tikzpicture}
    \node [mybox] (box){\usebox{\boiboite}};
    \node[fancytitle, right=10pt] at (box.north west) {\titre};
    \end{tikzpicture}
    \end{center}
    }

\newcommand{\newboxedtheorem}[4][]{%
    \couleurs[#1]
    \@ifnotempty{#4}{%
      \@ifundefined{the#4}{\@ifundefined{\thcounter}{\newcounter{#4}}{%
      \newcounter{#4}[\thcounter ] } } { }%
    }
    \newenvironment{#2}[1][]{%
    \@ifnotempty{#4}{\refstepcounter{#4}}
    \begin{boite}[#1]{\textbf{#3\@ifnotempty{#4}{ \csname the#4\endcsname}}\@ifnotempty{##1}{
    (##1)}}
    }%
    {%
    \end{boite}
    }
}
\makeatother
\newboxedtheorem{theorem}{Theorem}{theorem}
\newboxedtheorem{Th}{Theorem}{Th}          
\newboxedtheorem{Def}{Definition}{Def}     
\newboxedtheorem{Cor}{Corollary}{Cor}      
\newboxedtheorem{Pro}{Proposition}{Pro}    
\newboxedtheorem{Lem}{Lemma}{Lem}          
\newboxedtheorem{Exp}{Example}{Exp}        
\newboxedtheorem{Rem}{Remark}{Rem} 

\newboxedtheorem{SDem}{Suite de Preuve}{SDem}
\newboxedtheorem{SDe}{Suite de Preuve}{SDe}
\newboxedtheorem{Con}{Conséquence}{Con}
\newboxedtheorem{Cri}{Critère}{Cri}
\newboxedtheorem{test}{test}{test}
\newboxedtheorem{Prot}{Propriétés}{Prot}
\newboxedtheorem{Exer}{Exercice}{Exer}
\newboxedtheorem{suite}{Suite}{}

\newcommand{\N}{\mathbb{N}}
\newcommand{\R}{\mathbb{R}}

\usepackage{xcolor} 
\usepackage{hyperref}
\hypersetup{
    colorlinks=true,
    linkcolor=blue,      
    citecolor=blue,      
    urlcolor=blue }       

\begin{document}

\journaltitle{arXiv Preprint}
\copyrightyear{2025}
\pubyear{2025}
\date{20 April 2025}

\firstpage{1}

\title[Introduction to Hausdorff Measure and Dimension]{An Introduction to the Hausdorff Measure and Its Applications in Fractal Geometry}

\author[1,$\ast$]{Mohammed Nechba \ORCID{0009-0007-5850-3531}}
\author[2,$\ast$]{Mustapha Ouyaaz}
\author[1]{Abdellatif El Afia \ORCID{0000-0003-1921-4431}}
\author[2]{ Mohammed El Arrouchi \ORCID{0000-0002-7474-3697} }

\authormark{M. nechba et al.}

\address[1]{\orgdiv{École Nationale Supérieure d'Informatique et d'Analyse des Systèmes}, 
\orgname{Université Mohammed V de Rabat}, 
\orgaddress{\street{Avenue Mohamed Ben Abdellah, Av. Regragui}, \postcode{10000}, \state{Rabat}, \country{Morocco}}}

\address[2]{\orgdiv{Faculté des Sciences de Kénitra}, 
\orgname{Université Ibn Tofail}, 
\orgaddress{\street{B.P. 133, Avenue de l’Université}, \postcode{14000}, \state{Kénitra}, \country{Morocco}}}




\abstract{This paper presents a comprehensive introduction to the Hausdorff measure, a fundamental tool in fractal geometry and geometric measure theory. We begin by defining the Hausdorff outer measure on subsets of metric spaces, followed by a discussion of Carathéodory’s criterion, which characterizes measurable sets. From this foundation, we construct the Hausdorff measure and explore its essential properties, including monotonicity and translation invariance. We then introduce the Hausdorff dimension, a powerful generalization of Euclidean dimension, particularly suited to analyzing non-regular or self-similar sets. As an application, we examine the Cantor ternary set, computing its Hausdorff dimension and demonstrating how the Hausdorff measure captures its geometric complexity. This exposition aims to bridge the gap between abstract theory and illustrative application, offering insights relevant to mathematics and various scientific domains such as physics and computer science.}


\maketitle

\section{Introduction}
Fractal geometry has emerged as a fundamental framework for analyzing irregular structures in both pure mathematics and applied sciences \citep{Mandelbrot1982}. At the heart of this discipline lies the \emph{Hausdorff measure}, a powerful tool that extends classical geometric notions to highly irregular sets \citep{Falconer2014}. Building on the foundational work of \citet{Hausdorff1919}, this measure provides a rigorous approach to quantifying the "size" of fractal objects that defy traditional Euclidean analysis.

The Hausdorff measure $\mathcal{H}^s$ represents a profound generalization of conventional measures, characterized by its ability to:
\begin{itemize}
    \item Assign meaningful dimensions to non-smooth sets through the Hausdorff dimension $\dim_{\mathcal{H}}$
    \item Maintain essential geometric properties including:
    \begin{itemize}
        \item Translation invariance: $\mathcal{H}^s(A+x) = \mathcal{H}^s(A)$
        \item Scaling behavior: $\mathcal{H}^s(\lambda A) = \lambda^s\mathcal{H}^s(A)$
        \item Countable stability: $\dim_{\mathcal{H}}(\bigcup A_i) = \sup \dim_{\mathcal{H}}(A_i)$
    \end{itemize}
\end{itemize}

Our investigation begins with the construction of outer measures on metric spaces $(X,d)$, where we define:
$$\mu^*(A) = \inf\left\{\sum \text{diam}(U_i)^s : A \subset \bigcup U_i, \text{diam}(U_i) < \epsilon\right\}$$
and establish Carathéodory's crucial measurability criterion \citep{Rogers1998}:
$$ B \in \mathcal{M}(\mu^*) \iff \forall A \subset X, \mu^*(A) = \mu^*(A\cap B) + \mu^*(A\cap B^c) $$

The Hausdorff dimension emerges as a key invariant through the critical exponent:
$$\dim_{\mathcal{H}}(A) = \inf\{s \geq 0 : \mathcal{H}^s(A) = 0\} = \sup\{s \geq 0 : \mathcal{H}^s(A) = \infty\}$$

We demonstrate these concepts through a detailed analysis of the Cantor ternary set $C$, proving:
\begin{itemize}
    \item Its Lebesgue measure vanishes: $\lambda(C) = 0$
    \item Its Hausdorff dimension is exact: $\dim_{\mathcal{H}}(C) = \frac{\log 2}{\log 3}$
    \item The measure at critical dimension is non-trivial: $\frac{1}{2} \leq \mathcal{H}^s(C) \leq 1$ for $s = \dim_{\mathcal{H}}(C)$
\end{itemize}

Our treatment advances the existing literature by:
\begin{itemize}
    \item Presenting a unified development from outer measures to dimension theory
    \item Establishing sharp quantitative bounds for classical fractals
    \item Demonstrating the interplay between geometric and measure-theoretic properties
\end{itemize}

The paper's structure reflects this progression: Section \ref{sec:preliminaries} covers metric outer measures and measurability, Section \ref{sec:construction} develops the Hausdorff measure and dimension, and Section \ref{sec:applications} provides complete analysis of the Cantor set, including:
\begin{itemize}
    \item Its recursive construction $C = \bigcap_{n\geq 0} C_n$ where $C_n$ consists of $2^n$ intervals
    \item The exact computation $\mathcal{H}^s(C_n) = (\frac{2}{3^s})^n$
    \item The dimensional analysis via covering arguments
\end{itemize}

This work bridges the gap between abstract measure theory and concrete fractal analysis, offering insights valuable across mathematics, physics, and computational geometry \citep{Edgar2008,Falconer2013}.
\section{Preliminaries}
\label{sec:preliminaries}

We begin by recalling fundamental concepts from measure theory that will serve as the foundation for our construction of Hausdorff measure. These notions will allow us to generalize traditional geometric measurements to irregular fractal sets.

\subsection{Outer Measure}
\begin{Def}
Let $X$ be an arbitrary set. An \textbf{outer measure} on $X$ is a function $\mu^* : \mathcal{P}(X) \rightarrow [0, +\infty]$ satisfying:
\begin{itemize}
    \item i) $\mu^*(\varnothing) = 0$ (null empty set)
    \item ii) Monotonicity: If $A \subset B$, then $\mu^*(A) \leq \mu^*(B)$
    \item iii) Countable subadditivity: For any sequence $\{A_n\}_{n\in\mathbb{N}}$ of subsets of $X$,
    $$\mu^*\left(\bigcup_{n\in\mathbb{N}} A_n\right) \leq \sum_{n\in\mathbb{N}} \mu^*(A_n)$$
\end{itemize}
\end{Def}

\begin{remark}
An outer measure is not necessarily a positive measure unless it is additive. However, any positive measure defined on $\mathcal{P}(X)$ is automatically an outer measure. This distinction becomes crucial when working with non-measurable sets in the sense of Carathéodory.
\end{remark}

\begin{example}
Let $X$ be any set and define $\varphi : \mathcal{P}(X) \rightarrow \{0,1\}$ by:
$$
\varphi(A) = 
\begin{cases} 
0 & \text{if } A = \varnothing \\
1 & \text{otherwise}
\end{cases}
$$
Then $\varphi$ is an outer measure but not a positive measure. This simple example illustrates how outer measures can capture set existence rather than geometric size.
\end{example}

\subsection{Metric Outer Measure}
The abstract notion of outer measure gains geometric significance when we consider metric spaces. The additional structure allows us to relate set measurements to their spatial separation.

\begin{Def}
Let $(X,d)$ be a metric space. An outer measure $\mu^*$ is called a \textbf{metric outer measure} if for all $E, F \subset X$:
$$
d(E,F) > 0 \implies \mu^*(E \cup F) = \mu^*(E) + \mu^*(F)
$$
\vspace{-0.3cm}
where $d(E,F) = \inf\{d(x,y) : x \in E, y \in F\}$ is the separation distance between sets. This condition ensures the measure respects the metric topology.
\end{Def}

\begin{Def}
Let $X$ be a set equipped with an outer measure $\mu^*$. A subset $B \subset X$ is called \textbf{$\mu^*$-measurable} if for every subset $A \subset X$:
$$
\mu^*(A) = \mu^*(A \cap B) + \mu^*(A \cap B^\complement)
$$
The collection of all $\mu^*$-measurable sets is denoted by $\mathcal{M}(\mu^*)$. This measurability criterion, due to Carathéodory, will be essential for constructing proper measures from outer measures.
\end{Def} 

\subsection{Carathéodory's Theorem}
\begin{Th}[\cite{Caratheodory1914}, \cite{Folland1999}]
\label{th:1}
Let $\mu^*$ be an outer measure on $X$. Then:
\begin{itemize}

    \item The collection $\mathcal{M}(\mu^*)$ of all $\mu^*$-measurable subsets of $X$ forms a $\sigma$-algebra on $X$
    \item The restriction of $\mu^*$ to $\mathcal{M}(\mu^*)$ is a complete measure
\end{itemize}
\end{Th}

\vspace{-1.1cm}
\begin{proof}[Proof of Caratheodory’s Theorem \ref{th:1} : ]
We will denote $\mathcal{M} = \mathcal{M}(\mu^*)$ for simplicity.

$\bullet$ First, we show $\mathcal{M}$ is an algebra (closed under finite unions and complements, and contains $\varnothing$):

\begin{itemize}
    \item[-] For any $A \subset X$, since $\mu^*$ is an outer measure:
    \begin{align*}
    \mu^*(A \cap \varnothing) &= \mu^*(\varnothing) = 0 \\
    \mu^*(A \cap \varnothing^\complement) &= \mu^*(A \cap X) = \mu^*(A)
    \end{align*}
    Thus $\mu^*(A \cap \varnothing) + \mu^*(A \cap \varnothing^\complement) = \mu^*(A)$, proving $\varnothing \in \mathcal{M}$.
    
    \item[-] Let $B \in \mathcal{M}$ and $A \subset X$. Then:
    \begin{align*}
    \mu^*(A) &= \mu^*(A \cap B) + \mu^*(A \cap B^\complement) \\
             &= \mu^*(A \cap B^\complement) + \mu^*(A \cap (B^\complement)^\complement)
    \end{align*}
    Hence $B^\complement \in \mathcal{M}$.
    \item[-] Let $A_1, A_2 \in \mathcal{M}$ and $B \subset X$. Since $A_1 \in \mathcal{M}$:
    \vspace{-0.4cm}
    \begin{align*}
    \mu^*(B \cap (A_1 \cup A_2)) &= \mu^*(B \cap (A_1 \cup A_2) \cap A_1) + \mu^*(B \cap (A_1 \cup A_2) \cap A_1^\complement) \\
                                  &= \mu^*(B \cap A_1) + \mu^*(B \cap A_2 \cap A_1^\complement)
    \end{align*} 
    Therefore:
    \begin{align*}
    \mu^*(B \cap (A_1 \cup A_2)) &+ \mu^*(B \cap (A_1 \cup A_2)^\complement) \\
    &= \mu^*(B \cap A_1) + \mu^*(B \cap A_1^\complement \cap A_2) + \mu^*(B \cap A_1^\complement \cap A_2^\complement) \\
    &= \mu^*(B \cap A_1) + \mu^*(B \cap A_1^\complement) \quad (\text{since } A_2 \in \mathcal{M}) \\
    &= \mu^*(B) \quad (\text{since } A_1 \in \mathcal{M})
    \end{align*}
    Thus $A_1 \cup A_2 \in \mathcal{M}$.
\end{itemize}

$\bullet$ To show $\mathcal{M}$ is a $\sigma$-algebra, it suffices to prove closure under countable unions of pairwise disjoint sets:

Let $(B_n)_{n \geq 1}$ be pairwise disjoint sets in $\mathcal{M}$. We prove by induction that for all $A \subset X$ and $n \geq 1$:
$$\mu^*(A) = \sum_{k=1}^n \mu^*(A \cap B_k) + \mu^*\left(A \setminus \bigcup_{k=1}^n B_k\right)$$

\begin{itemize}
    \item[-] Base case ($n=1$): Follows directly from $B_1 \in \mathcal{M}$
    \item[-] Inductive step: Assume holds for $n$, then for $B_{n+1} \in \mathcal{M}$:
    \begin{align*}
    \mu^*\left(A \setminus \bigcup_{k=1}^n B_k\right) &= \mu^*\left(\left(A \setminus \bigcup_{k=1}^n B_k\right) \cap B_{n+1}\right) \\
    &\quad + \mu^*\left(\left(A \setminus \bigcup_{k=1}^n B_k\right) \setminus B_{n+1}\right) \\
    &= \mu^*(A \cap B_{n+1}) + \mu^*\left(A \setminus \bigcup_{k=1}^{n+1} B_k\right)
    \end{align*}
    Combining with the inductive hypothesis yields the result for $n+1$
\end{itemize}

Now let $B = \bigcup_{n \geq 1} B_n$. We have:
\begin{align*}
\mu^*(A) &\geq \sum_{n \geq 1} \mu^*(A \cap B_n) + \mu^*(A \setminus B) \quad \text{(by monotonicity)} \\
&\geq \mu^*\left(\bigcup_{n \geq 1} (A \cap B_n)\right) + \mu^*(A \setminus B) \quad \text{(subadditivity)} \\
&= \mu^*(A \cap B) + \mu^*(A \setminus B)
\end{align*}

The reverse inequality follows from subadditivity, so equality holds, proving $B \in \mathcal{M}$.

For arbitrary countable unions, decompose into disjoint sets as usual.

$\bullet$ Finally, $\mu^*$ restricted to $\mathcal{M}$ is a measure:
\begin{itemize}
    \item[-] $\mu^*(\varnothing) = 0$ by definition
    \item[-] For disjoint $(B_n)_{n \geq 1}$ in $\mathcal{M}$, take $A = \bigcup_{n \geq 1} B_n$ in the earlier equality:
    $$\mu^*\left(\bigcup_{n \geq 1} B_n\right) = \sum_{n \geq 1} \mu^*(B_n)$$
    proving countable additivity.
\end{itemize}
\end{proof}
\begin{Pro}
\label{pro:1}
Let $(X, d)$ be a metric space, and $\mu^*$ an outer measure on $X$.\\
If $\mu^*$ is a metric outer measure, then $\mathcal{B}(X) \subset \mathcal{M}(\mu^*)$, 
where $\mathcal{B}(X)$ is the Borel $\sigma$-algebra.
\end{Pro}

The proof of this proposition is based on the following lemma:

\begin{Lem}
\label{lem:1}
Let $\mu^*$ be a metric outer measure, and $\{A_n\}_{n\in\mathbb{N}}$ 
an increasing sequence of subsets of $X$. Let $A =\bigcup_{n\in\mathbb{N}} A_n$. 
If $d(A_n, A \setminus A_{n+1}) > 0$ for all $n\in\mathbb{N}$, then we have
$$\mu^*(A) = \lim_{n \to +\infty} \mu^*(A_n)$$
\end{Lem}

\begin{proof}[Proof of Lemma \ref{lem:1}: ]

Let $\mu^*$ be a metric outer measure, and $\{A_n\}_{n\in\mathbb{N}}$ 
an increasing sequence of subsets of $X$. Let $A =\bigcup_{n\in\mathbb{N}} A_n$. Assume $d(A_n, A \setminus A_{n+1}) > 0$ for all $n\in\mathbb{N}$. For $\mu^*(A) < +\infty$.\\ 
Let $B_0 = A_0$ and for $n\geq 1$, $B_n = A_n \setminus A_{n-1}$. 
Since $\{A_n\}_{n\in\mathbb{N}}$ is increasing, $\{B_n\}_{n\in\mathbb{N}}$ 
is a pairwise disjoint family of subsets of $X$ satisfying
$$\bigcup_{n\in\mathbb{N}} B_n =\bigcup_{n\in\mathbb{N}} A_n =A$$ 

For all $n\in\mathbb{N}$, we have
$$A_n \subset A \subset A_n\cup\bigcup_{k\geq n+1}B_k$$ 
Therefore, by monotonicity and subadditivity of $\mu^*$:
\begin{align}
\mu^*(A_n) \leq \mu^*(A) \leq \mu^*(A_n)+\sum_{k\geq n+1} \mu^*(B_k)
\end{align}

Thus, to show $\lim_{n \to +\infty} \mu^*(A_n)=\mu^*(A)$, 
it suffices to show $\lim_{n \to +\infty} \sum_{k\geq n+1}\mu^*(B_k)=0$. 
This reduces to proving the series $\sum \mu^*(B_k)$ converges.

Indeed, suppose there exist $n$ and $m$ of the same parity with $n < m$ 
and $d(B_n, B_m)=0$. Then:
$$\forall \varepsilon >0,\;\exists (x,y) \in B_n \times B_m : d(x, y) <\varepsilon$$
Take $\varepsilon = d(A_n, A \setminus A_{n+1}) >0$ by hypothesis. 
There exists $(x_0, y_0) \in B_n\times B_m$ with $d(x_0, y_0) <\varepsilon$.

Now, $x_0 \in B_n \subset A_n$ and $y_0\in B_m=A_m \setminus A_{m-1}$.\\ 
Since $n+1 \leq m - 1$ (same parity and $n < m$) and $\{A_n\}_{n\in\mathbb{N}}$ is increasing, 
then $y_0\in A_m\setminus A_{m-1} \subset A_m\setminus A_{n+1} \subset A\setminus A_{n+1}$, 
hence $y_0\in A\setminus A_{n+1}$.\\ 
Thus $d(x_0, y_0) \geq d(A_n, A \setminus A_{n+1})= \varepsilon$. Contradiction.\\ 

Therefore, by hypothesis, for all $N\in \mathbb{N}$:
$$\mu^*\left(\bigcup_{k=0}^N B_{2k}\right) =\sum_{k=0}^N \mu^*(B_{2k}),\;
\mu^*\left(\bigcup_{k=0}^N B_{2k+1}\right) =\sum_{k=0}^N \mu^*(B_{2k+1})$$

By monotonicity of $\mu^*$:
$$\mu^*\left(\bigcup_{k=0}^N B_{2k}\right)\leq \mu^*(A) <\infty$$
and
$$\mu^*\left(\bigcup_{k=0}^N B_{2k+1}\right)\leq \mu^*(A) <\infty$$

Thus for all $N\in \mathbb{N}$:
\begin{align*}
\sum_{k=0}^{2N+1} \mu^*(B_{k}) &=\sum_{k=0}^N \mu^*(B_{2k})+\sum_{k=0}^N \mu^*(B_{2k+1})\\
&=\mu^*\left(\bigcup_{k=0}^N B_{2k}\right)+\mu^*\left(\bigcup_{k=0}^N B_{2k+1}\right)\\
&\leq 2\mu^*(A)
\end{align*}

Taking limit as $N \to \infty$:
$$\sum_{k\geq0} \mu^*(B_{k})\leq 2\mu^*(A)<\infty$$

Finally:
$$\mu^*(A) = \lim_{n \to +\infty} \mu^*(A_n)$$
\end{proof}

\begin{proof}[Proof of Proposition \ref{pro:1}: ]
Let $(X, d)$ be a metric space and $\mu^*$ an outer measure on $X$.\\ 
Assume $\mu^*$ is metric, i.e., $\forall E,F \subset X$ with $d(E,F) > 0$, 
we have $\mu^*(E\cup F)= \mu^*(E)+ \mu^*(F)$.\\ 

Let $F$ be a closed set in $X$. 
We show $F$ is $\mu^*$-measurable. For any subset $A \subset X$, 
we want to show $\mu^*(A)\geq \mu^*(A\cap F)+\mu^*(A\cap F^\complement)$, 
the reverse inequality following from subadditivity. 
We may assume $\mu^*(A)<\infty$, otherwise trivial.

For all $n\geq 1$, define:
$$F_n=\{x\in X \mid d(x,F)> \frac{1}{n}\}$$

Since $d(A \cap F, A \cap F_n ) > 0$, we have:
$$\mu^*((A \cap F_n)\cup (A \cap F))= \mu^*(A \cap F_n)+ \mu^*(A \cap F)$$

By monotonicity of $\mu^*$:
\begin{align*}
\mu^*(A) &\geq\mu^*((A \cap F_n)\cup (A \cap F))\\
         &= \mu^*(A \cap F_n)+ \mu^*(A \cap F)
\end{align*}

Taking limit as $n\to\infty$, since $\bigcup_{n\geq 1} F_n =F^\complement$, 
and $\{A\cap F_n\}_{n\geq 1}$ satisfies the conditions of the previous lemma, 
we obtain the desired inequality.\

Finally, all closed sets are $\mu^*$-measurable, 
and since closed sets generate the Borel $\sigma$-algebra, 
we have $\mathcal{B}(X) \subset \mathcal{M}(\mu^*)$.
\end{proof}

\subsection{Construction of Measures}
In this section, we focus on the construction of outer measures. For this purpose, let E be a set, $\mathcal{U} \subset \mathcal{P}(E)$ a covering of E, and $\delta : \mathcal{U} \longrightarrow [0, +\infty]$ a function with $\delta(\emptyset) = 0$. We define:
$$\mu(A) = \inf_{D \in \mathcal{R}(A)} \sum_{X \in D} \delta(X),\ \forall A \subset E$$
where $\mathcal{R}(A)$ denotes the set of countable coverings of A by elements of $\mathcal{U}$.

\begin{Th}
\label{th:2}
The mapping $\mu$ is an outer measure on E.
\end{Th}
\vspace{-0.4cm}
\begin{proof}[Proof of Theorem \ref{th:2}:]
\emph{1. Non-negativity:} For any $A \subset E$, by construction of $\mu$:
$$\mu(A) = \inf_{D \in \mathcal{R}(A)} \sum_{X \in D} \delta(X)$$
Since $\delta(X) \geq 0$ for all $X \in D$ and $D \in \mathcal{R}(A)$, we have:
$$\forall D \in \mathcal{R}(A),\ \sum_{X \in D} \delta(X) \geq 0$$
Therefore:
$$\mu(A) = \inf_{D \in \mathcal{R}(A)} \sum_{X \in D} \delta(X) \geq 0$$
Hence $\mu$ is well-defined.\\
\vspace{-0.3cm}
\emph{2. Null empty set:}  Verify that $\mu(\emptyset) = 0$: we have,
$$\mu(\emptyset) = \inf_{D \in \mathcal{R}(\emptyset)} \sum_{X \in D} \delta(X) \geq 0$$

Since $\{\emptyset\}$ is a covering of $\emptyset$, then $\{\emptyset\} \in \mathcal{R}(\emptyset)$. Thus:
\begin{align*}
\mu(\emptyset) &\leq \sum_{X \in \{\emptyset\}} \delta(X) \\
               &= \delta(\emptyset) \\
               &= 0
\end{align*}

Hence $\mu(\emptyset) = 0$.

\emph{3. Monotonicity:}  Let $A,B \subset E$ with $A \subset B$ and $\{A_n\}_{n\in\mathbb{N}} \in \mathcal{R}(B)$. Since $B \subset \bigcup_{n\in\mathbb{N}} A_n$ and $A \subset B$, we have $A \subset \bigcup_{n\in\mathbb{N}} A_n$. Therefore:
\vspace{-0.5cm}
$$\{A_n\}_{n\in\mathbb{N}} \in \mathcal{R}(A)$$

Consequently:
\vspace{-0.5cm}
$$\mathcal{R}(B) \subset \mathcal{R}(A)$$

Which implies:
\vspace{-0.6cm}
$$\inf_{D \in \mathcal{R}(A)} \sum_{X \in D} \delta(X) \leq \inf_{D \in \mathcal{R}(B)} \sum_{X \in D} \delta(X)$$

Thus $\mu(A) \leq \mu(B)$.

\emph{4. Countable subadditivity:} Let $\{A_n\}_{n\in\mathbb{N}}$ be a family of subsets of E and $\epsilon > 0$. For each n, choose $D_n \in \mathcal{R}(A_n)$ such that:
\vspace{-0.5cm}
$$\sum_{X \in D_n} \delta(X) \leq \mu(A_n) + \frac{\epsilon}{2^{n+1}}$$

The existence of such coverings is guaranteed by the construction of $\mu$. Since $\bigcup_{n\in\mathbb{N}} D_n \in \mathcal{R}(\bigcup_{n\in\mathbb{N}} A_n)$, we have:
\vspace{-0.5cm}
\begin{align*}
\mu\left(\bigcup_{n\in\mathbb{N}} A_n\right) &\leq \sum_{X \in \bigcup_{n\in\mathbb{N}} D_n} \delta(X) \\
                                  &\leq \sum_{n\in\mathbb{N}} \sum_{X \in D_n} \delta(X) \\
                                  &\leq \sum_{n\in\mathbb{N}} \left(\mu(A_n) + \frac{\epsilon}{2^{n+1}}\right) \\
                                  &= \epsilon + \sum_{n\in\mathbb{N}} \mu(A_n)
\end{align*}
\vspace{-0.5cm}
Taking $\epsilon \to 0$ establishes the subadditivity of $\mu$. Therefore $\mu$ is an outer measure.
\end{proof}

\section{Hausdorff Measure}
\label{sec:construction}
\subsection{Definition of Hausdorff Measure}

We carry out here the same measure construction as before, working on a metric space (E, d), which we cover with the class $\mathcal{U}=\{X\subset E ,\text{diam}(X)\leq\epsilon\}$, where $\epsilon$ is strictly positive. We take a real number $s > 0$ and define the mapping $\delta: X\longrightarrow \text{diam}(X)^s$.

Finally, we define $\mathcal{R}_{\epsilon}(A)$, the set of countable coverings of A by subsets $X\subset E$ such that $\text{diam}(X)\leq \epsilon$. We then obtain the quantity:

$$\forall A \subset E,\mathcal{H}_{\epsilon}^s(A):= \inf_{D \in \mathcal{R}_{\epsilon}(A)} \sum_{X\in D} (\text{diam}(X))^s$$

\begin{Def}
We define the Hausdorff measure as follows:

$$\mathcal{H}^s(A):=\lim_{\epsilon\to 0} \mathcal{H}_{\epsilon}^s(A),\ \forall A \subset E$$
\end{Def}

\begin{remark}

\begin{itemize}
    \item[-] Why does the limit exist? \\
    The mapping $\eta >0 \rightarrow \mathcal{H}_{\eta}^s(A)$, for all A $\subset$ X, is decreasing (since the infimum is taken over coverings with smaller and smaller diameters) and therefore admits a limit at 0.
\item[-] Note that if $\eta <\epsilon$, then $\mathcal{R}_{\eta}(A)\subset \mathcal{R}_{\epsilon}(A)$, so $\mathcal{H}_{\eta}^s(A)\geq \mathcal{H}_{\epsilon}^s(A)$. \\Consequently, $\mathcal{H}^s(A)=\sup_{\epsilon >0}\mathcal{H}_{\epsilon}^s(A)$ is well-defined.
\end{itemize}

\end{remark}

\begin{Pro}
\label{pro:2}
The Hausdorff measure is:
\begin{itemize}
\item i) an outer measure.
\item ii) a metric measure.
\item iii) a Borel measure.
\end{itemize}
\end{Pro}

\begin{proof}[Proof of Proposition \ref{pro:2}: ]
\textbf{ i)} By construction, $\mathcal{H}^s$ is an outer measure.\\
\textbf{ ii)} Let $s\geq 0$, $A,B\subset E$ with $d(A,B)>0$, and $\epsilon>0$ with $d(A,B)>\epsilon$.
Consider a covering $D \in \mathcal{R}_{\epsilon}(A\cup B)$. If $X \in D$, the set X cannot intersect both A and B simultaneously. Consequently, $D=D_1 \cup D_2$, $D_1 \cap D_2 =\varnothing$ and $D_1 \in \mathcal{R}_{\epsilon}(A)$, $D_2 \in \mathcal{R}_{\epsilon}(B)$. We deduce that
\vspace{-0.5cm}
$$\mathcal{H}^s(A)+ \mathcal{H}^s(B) \leq \sum_{X\in D_1} (\text{diam}(X))^s+\sum_{X\in D_2} (\text{diam}(X))^s=\sum_{X\in D} (\text{diam}(X))^s$$
\vspace{-0.5cm}
Taking the infimum over $D \in \mathcal{R}_{\epsilon}(A\cup B)$ and letting $\epsilon$ tend to 0, we obtain:
\vspace{-0.3cm}
$$\mathcal{H}^s(A)+\mathcal{H}^s(B) \leq \mathcal{H}^s(A\cup B)$$
\vspace{-0.5cm}
By subadditivity of $\mathcal{H}^s$ we have:
\vspace{-0.5cm}
$$\mathcal{H}^s(A\cup B) \leq \mathcal{H}^s(A)+\mathcal{H}^s(B)$$

Consequently,
\vspace{-0.5cm}
$$\mathcal{H}^s(A\cup B) =\mathcal{H}^s(A)+\mathcal{H}^s(B)$$

It follows that $\mathcal{H}^s$ is indeed a metric measure.\\
\textbf{ iii)} By Carathéodory's Theorem, $\mathcal{H}^s$ is a measure on the $\sigma$-algebra $M(\mathcal{H}^s)$, and since we have $\mathcal{B}(E)\subset M(\mathcal{H}^s)$ (proposition \ref{pro:1}). Therefore $\mathcal{H}^s$ is a Borel measure.
\end{proof}

\subsection{Some Properties of Hausdorff Measure}
\begin{Pro}
\label{pro:3}
\begin{itemize}
\item i) The Hausdorff measure is translation invariant.
\item ii) The Hausdorff measure satisfies the following scaling property:
\end{itemize}
$$ \mathcal{H}^s(\lambda A) =\lambda^s \mathcal{H}^s(A),\;\;\;\;\;\;\;\;\forall s\in [0,+\infty]\;,\forall A\subset E\; ,\forall \lambda>0$$
\end{Pro}
\vspace{-1cm}
\begin{proof}[Proof of Proposition \ref{pro:3}: ]
Let $s \geq 0$, $\epsilon > 0$ and $A \subset E$.

\textbf{i)} Let $\{A_n\}_{n\in \mathbb{N}}\in \mathcal{R}_{\epsilon}(A)$ then for all $x \in E$, $\{A_n +x\}_{n\in \mathbb{N}}$ is a covering of $A+x$ where
$A+x=\{a+x ,a\in A\}$, moreover for all $n\in \mathbb{N}$ we have:
\vspace{-0.6cm}
\begin{align*}
\text{diam}(A_n +x)=\text{diam}(A_n) &\Rightarrow \text{diam}(A_n +x)^s=\text{diam}(A_n)^s\\
&\Rightarrow \sum_{n\in \mathbb{N}}\text{diam}(A_n +x)^s= \sum_{n\in \mathbb{N}}\text{diam}(A_n)^s\\
& \Rightarrow \inf_{\{A_n\}_{n\in \mathbb{N}}\in \mathcal{R}_{\epsilon}(A)} \sum_{n\in \mathbb{N}}\text{diam}(A_n +x)^s= \inf_{\{A_n\}_{n\in \mathbb{N}}\in \mathcal{R}_{\epsilon}(A)}\sum_{n\in \mathbb{N}}\text{diam}(A_n)^s\\
&\Rightarrow \mathcal{H}_{\epsilon}^s(A+x)=\mathcal{H}_{\epsilon}^s(A)
\end{align*}
\vspace{-0.3cm}
Taking the limit as $\epsilon$ tends to zero we obtain $\mathcal{H}^s(A+x)=\mathcal{H}^s(A)$.

\textbf{ii)} Let $\lambda >0$. For any covering $\{A_n\}_{n\in \mathbb{N}}$ of A where $\text{diam}(A_n)\leq\epsilon$ corresponds a covering $\{\lambda A_n\}_{n\in \mathbb{N}}$ of $\lambda A$, where $\lambda A=\{\lambda a ,a\in A\}$ and $\text{diam}(\lambda A_n)\leq \lambda \epsilon$. We have 
\begin{align*}
\text{diam}(\lambda A_n)= \lambda \text{diam}(A_n) &\Rightarrow \text{diam}(\lambda A_n)^s= \lambda^s \text{diam}(A_n)^s\\
&\Rightarrow \sum_{n\in \mathbb{N}} \text{diam}(\lambda A_n)^s= \sum_{n\in \mathbb{N}} \lambda^s \text{diam}(A_n)^s
\end{align*}
Taking the infimum over all coverings $\{A_n\}_{n\in \mathbb{N}}$ we obtain $\mathcal{H}_{\lambda \epsilon}^s(\lambda A)=\lambda^s \mathcal{H}_{\epsilon}^s(A)$, and taking $\epsilon$ to zero gives
$\mathcal{H}^s(\lambda A)=\lambda^s \mathcal{H}^s(A)$

\end{proof}

\begin{Pro}
\label{pro:4}
\begin{itemize}
\item i) The measure $\mathcal{H}^0$ is the counting measure.
\item ii) On $\mathbb{R}$, the measure $\mathcal{H}^1$ is the Lebesgue measure.
\end{itemize}
\end{Pro}
\begin{proof}[Proof of Proposition \ref{pro:4}: ]
 \textbf{i)} Let A be a finite subset of E and let $\delta_0 =\min \{d(x,y) \diagup x,y\in A ,x \neq y \}$. We want to show that for all $\delta \in ]0, \delta_0[$, $\mathcal{H}_{\delta}^0(A) =\text{card}(A)$. Indeed, since $A \subset \bigcup_{x\in A} \{x\}$, then:
$$\mathcal{H}_{\delta}^0(A)\leq \sum_{x\in A} 1 =\text{card}(A)$$

Conversely, for all $\delta < \delta_0$, if $F \in\mathcal{R}_{\delta}(A)$, any set in F contains at most one element of A, so F consists of at least n subsets of E with $n=\text{card}(A)$ and thus $\mathcal{H}_{\delta}^0(A) \geq \text{card}(A)$. This shows that for all sufficiently small $\delta$, $\mathcal{H}_{\delta}^0(A) = \text{card}(A)$ and therefore $\mathcal{H}^0(A)= \text{card}(A)$.

 If A is an infinite subset of E, it contains finite sets with arbitrarily many elements. By monotonicity of $\mathcal{H}^0$, we have $\mathcal{H}^0(A) \geq n$ for all $n \in \mathbb{N}$ and thus $\mathcal{H}^0(A) =+\infty$. This is indeed the counting measure on $E$.

\textbf{ii)} Since the Hausdorff measure is translation invariant, it suffices to show that $\mathcal{H}^1([0,1])=1$. First, we show that $\mathcal{H}^1([0,1])\leq 1$.

Let $U_i =[\frac{i-1}{k} ;\frac{i}{k}]$ for $i=1,...,k$. Then $[0,1]=\bigcup_{i=1}^k U_i$ and $\text{diam}(U_i)=\frac{1}{k}$ for $i=1,...,k$, so:
$$ \mathcal{H}_{\frac{1}{k}}^1([0,1])\leq \sum_{i=1}^k \text{diam}(U_i)=1 $$

Taking k to infinity gives $\mathcal{H}^1([0,1])\leq 1$.

For the reverse inequality, consider $\epsilon >0$, $\delta>0$ and $\{A_n\}_{n\in \mathbb{N}} \in \mathcal{R}_{\delta}([0,1])$ such that:
\vspace{-0.5cm}
$$\mid \mathcal{H}_{\delta}([0,1]) -\sum_{n\in \mathbb{N}} \text{diam}(A_n) \mid \leq \epsilon $$

Without loss of generality, we may assume the $A_n$ are open intervals. There exists a family $\{a_k\}_{0\leq k\leq n}$ with $0=a_0 < a_1<...<a_n=1$ such that $[a_k,a_{k+1}]$ is contained in exactly one $A_k$ for all k=0,...,n-1. Consequently:
\vspace{-0.5cm}
$$ \mathcal{H}_{\delta}^1([0,1]) +\epsilon \geq \sum_{k=0}^{n-1} \text{diam}([a_k,a_{k+1}])= 1$$
Taking $\epsilon$ and $\delta$ to 0 gives $\mathcal{H}^1([0,1]) \geq 1$. It follows that $\mathcal{H}^1([0,1]) = 1$.
\end{proof}

\begin{Pro}
\label{pro:5}
Let f be a similarity with ratio $r > 0$, that is:\\ $\forall x,y \in E, d(f(x),f(y))=rd(x,y)$. Then:
\vspace{-0.5cm}
$$ \forall s\geq 0,\; \forall A\subset E,\mathcal{H}^s(f(A))=r^s\mathcal{H}^s(A)$$
\end{Pro}

\begin{proof}[Proof of Proposition \ref{pro:5}: ]
Let $A \subset E$, $\epsilon >0$. Note that $f(\mathcal{R}_{\epsilon}(A)) = \mathcal{R}_{\epsilon r}(f(A))$. We deduce:
\vspace{-0.5cm}
\begin{align*}
\mathcal{H}_{\epsilon r}^s(f(A))&=\inf_{D \in f(\mathcal{R}_{\epsilon}(A))} \sum_{X\in D} (\text{diam}(X))^s\\
&=\inf_{D' \in \mathcal{R}_{\epsilon}(A)} \sum_{X\in D'} (\text{diam}(f(X)))^s\\
&=r^s \inf_{D' \in \mathcal{R}_{\epsilon}(A)} \sum_{X\in D'} (\text{diam}(X))^s\\
&=r^s \mathcal{H}_{\epsilon}^s(A)
\end{align*}

Taking $\epsilon$ to 0 gives $\mathcal{H}^s(f(A))=r^s \mathcal{H}^s(A)$.
\end{proof}

\begin{Pro}
\label{pro:6}
Let $(X, d_x)$ and $(Y, d_y)$ be two metric spaces, $k > 0$ and $f : (X, d_x)\longrightarrow (Y, d_y)$ a k-Lipschitz function.Then:
$$\forall s\geq 0,\; \forall A\subset X,\mathcal{H}^s(f(A))\leq k^s\mathcal{H}^s(A)$$
\end{Pro}
\vspace{-1cm}
\begin{proof}[Proof of Proposition \ref{pro:6}: ]
By the same proof method as above, we have:
\vspace{-0.5cm}
\begin{align*}
\mathcal{H}_{\epsilon}^s(f(A))&=\inf_{D \in f(\mathcal{R}_{\epsilon/k}(A))} \sum_{X\in D} (\text{diam}(X))^s\\
&=\inf_{D' \in \mathcal{R}_{\epsilon/k}(A)} \sum_{X\in D'} (\text{diam}(f(X)))^s\\
&\leq k^s \inf_{D' \in \mathcal{R}_{\epsilon/k}(A)} \sum_{X\in D'} (\text{diam}(X))^s\\
&\leq k^s \mathcal{H}_{\epsilon/k}^s(A)
\end{align*}
Taking $\epsilon$ to 0 gives $\mathcal{H}^s(f(A))\leq k^s \mathcal{H}^s(A)$.
\end{proof}

\subsection{Definition of Hausdorff Dimension}
\begin{Pro}
\label{pro:7}
\begin{itemize}
\item i) Let $A\subset \R^n$. If there exists $d \geq 0$ such that $\mathcal{H}^d(A)<\infty$, then:
\vspace{-0.5cm}
$$\mathcal{H}^s(A)=0,\;\;\;\; \forall s>d$$
\item ii) For any subset $A\subset \R^n$, there exists a unique value $d \in \mathbb{R}^+ \cup \{\infty\}$ such that:
\vspace{-0.5cm}
$$\mathcal{H}^s(A)=+\infty,\;\;\;\; \forall s<d$$
and
\vspace{-0.5cm}
$$\mathcal{H}^s(A)=0,\;\;\;\;\;\; \forall s>d$$
\end{itemize}
\end{Pro}

\begin{proof}[Proof of Proposition \ref{pro:7}: ]
\textbf{i)} Let $A\subset \R^n$, $\epsilon > 0$ and $s>d\geq 0$:
\begin{align*}
\mathcal{H}_{\epsilon}^s(A) &= \inf_{D \in \mathcal{R}_{\epsilon}(A)} \sum_{X\in D} (\text{diam}(X))^s\\
&= \inf_{D \in \mathcal{R}_{\epsilon}(A)} \sum_{X\in D}(\text{diam}(X))^d(\text{diam}(X))^{s-d}\\
&\leq \inf_{D \in \mathcal{R}_{\epsilon}(A)} \sum_{X\in D}(\text{diam}(X))^d \epsilon^{s-d}\\
&= \epsilon^{s-d} \mathcal{H}_{\epsilon}^d(A)
\end{align*}
Since $\mathcal{H}^d(A)<\infty$, we obtain $\mathcal{H}_{\epsilon}^s(A)\to 0$ as $\epsilon \to 0$, and thus $\mathcal{H}^s(A)=0$ for all $A\subset E$ when $s>d$.\\

\textbf{ii)} Existence:\\
If $\{\mathcal{H}^s(A)<+\infty\} = \emptyset$, then we set $d=+\infty$, which satisfies the theorem's conditions. Otherwise, we set $d=\inf\{s\geq 0 : \mathcal{H}^s(A)<+\infty\}$.\\
If $s<d$ then $\mathcal{H}^s(A)=+\infty$.\\
If $s>d$, then there exists $s_0 \in ]d,s[$ such that $\mathcal{H}^{s_0}(A)<+\infty$. By part \textbf{i)}, we deduce that $\mathcal{H}^s(A)=0$, which proves existence.Uniqueness:
By contradiction, if there existed $s_1 > s_2$ satisfying both hypotheses, then for $t \in ]s_1, s_2[$, the t-dimensional measure of $\R^n$ would be both $0$ and $\infty$.
\end{proof}
\begin{Def}
For any $A\subset \R^n$, the Hausdorff dimension of A is the value:
\vspace{-0.5cm}
\begin{align*}
\dim_{\mathcal{H}}(A) &= \sup\{s\geq 0 : \mathcal{H}^s(A)=+\infty\}\\
&= \inf\{s\geq 0 : \mathcal{H}^s(A)=0\}
\end{align*}
\end{Def}

\begin{figure}[H]
    \centering
    \includegraphics[scale=0.3]{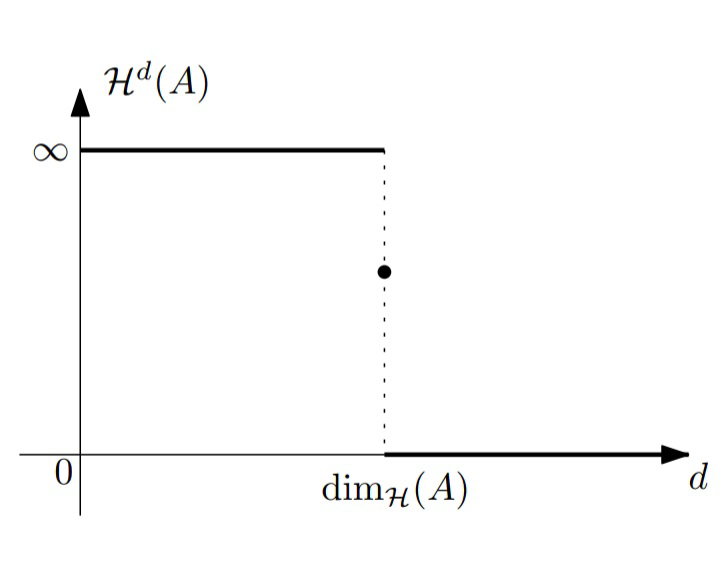}
    \caption{\textbf{Graph of the function }$d\rightarrow \mathcal{H}^d(A)$}
    \label{fig:enter-label}
\end{figure}

\subsection{Properties of Hausdorff Dimension}
\begin{Pro}
\label{pro:8}
\begin{itemize}
\item i) Let $A,B\subset \mathbb{R}^n$. If $A\subset B$ then $\dim_{\mathcal{H}}(A)\leq \dim_{\mathcal{H}}(B)$
\item ii) Let $f : \R^n \longrightarrow \R^n$ be a similarity, then $\dim_{\mathcal{H}}(f(A)) = \dim_{\mathcal{H}}(A)$
\end{itemize}
\end{Pro}

\begin{proof}[Proof of Proposition 8; ]
 Let A,B be subsets of $\R^n$ such that $A\subset B$ and $d\geqslant 0$ with $\mathcal{H}^d(B)=0$.\\
Since $\mathcal{H}^s$ is an outer measure on $\R^n$, we have $\mathcal{H}^s(A)\leqslant \mathcal{H}^s(B)$ for all $s\geqslant 0$.\\
In particular $\mathcal{H}^d(A)\leqslant \mathcal{H}^d(B)=0$, hence $\mathcal{H}^d(A)=0$.\\
Therefore, if $d\in \{s\geqslant 0 : \mathcal{H}^s(B)=0\}$ then $d\in \{s\geqslant 0 : \mathcal{H}^s(A)=0\}$. Thus,
$$\{s\geqslant 0 : \mathcal{H}^s(B)=0\} \subset \{s\geqslant 0 : \mathcal{H}^s(A)=0\}$$
It follows that:
$$\inf\{s\geqslant 0 : \mathcal{H}^s(A)=0\} \leqslant \inf\{s\geqslant 0 : \mathcal{H}^s(B)=0\}$$
By definition of Hausdorff dimension, we conclude that $\dim_{\mathcal{H}}(A)\leq \dim_{\mathcal{H}}(B)$.
\vspace{-1cm}
\begin{align*}
\dim_{\mathcal{H}}(f(A)) &= \inf\{s\geqslant 0 : \mathcal{H}^s(f(A))=0\} \;\;\;\;\;\;\; \text{(by definition)}\\
&= \inf\{s\geqslant 0 : r^s\mathcal{H}^s(A)=0\} \;\;\;\;\;\;\; \text{(Proposition 5)}\\
&= \inf\{s\geqslant 0 : \mathcal{H}^s(A)=0\}\\
&= \dim_{\mathcal{H}}(A).
\end{align*}
\end{proof}

\begin{Th}
\label{th:3}
Any countable subset A of $\R^n$ has Hausdorff dimension zero.
\end{Th}
\vspace{-0.5cm}
\begin{proof}[Proof of Theorem \ref{th:3}: ]
It suffices to prove the theorem for countable subsets of $\R^n$.\\
Let $A=\{x_k\in\R^n : k\in\N\}$. We have $A\subset\bigcup_{k\in \N}\{x_k\}$ and $\text{diam}(\{x_k\})=0\leq \varepsilon$ for all $k\in \N$ and $\varepsilon>0$. This shows $(\{x_k\})_{k\in \N}\in \mathcal{R}_\varepsilon(A)$.

Consequently, for all $s>0$ and $\varepsilon >0$, we have
\vspace{-0.5cm}
$$\mathcal{H}_{\varepsilon}^s(A)\leq \sum_{k\in \N} (\text{diam}(\{x_k\}))^s =0$$
Taking the limit as $\varepsilon\rightarrow 0$, we obtain
\vspace{-0.5cm}
$$\mathcal{H}^s(A)=\lim_{\varepsilon\rightarrow 0}\mathcal{H}_{\varepsilon}^s(A)=0\;\;\;\;\; \forall s>0$$
It follows that $\dim_{\mathcal{H}}(A)=0$ (by definition of Hausdorff dimension).
\end{proof}

\begin{Pro}
\label{pro:9}
Let $(A_k)_{k\in \N}$ be a sequence of subsets of $\R^n$. Then:
\vspace{-0.5cm}
$$\dim_{\mathcal{H}}(\bigcup_k A_k)=\sup \{\dim_{\mathcal{H}}(A_k):k\in \N\}$$
\end{Pro}

\begin{proof}[Proof of Proposition \ref{pro:9}: ]
Let $s=\sup\{\dim_{\mathcal{H}}(A_n) :n\in\N\}$. We consider two cases:

\textit{First case:} If $t>s$, then for all $n\in \N$, $t>\dim_{\mathcal{H}}(A_n)$.

This implies $\mathcal{H}^t(A_n)=0$, and by subadditivity of $\mathcal{H}^t$:
$$\mathcal{H}^t(\bigcup_{n\in \N}A_n)\leq \sum_{n\in \N} \mathcal{H}^t(A_n)=0$$
Hence $\mathcal{H}^t(\bigcup_{n\in \N}A_n)=0$.

\textit{Second case:} If $t<s$, then by definition of supremum there exists $n_0 \in \N$ such that $\dim_{\mathcal{H}}(A_{n_0}) >t$, so $\mathcal{H}^t(A_{n_0}) =+\infty$.\\
By monotonicity of $\mathcal{H}^t$, since $A_{n_0}\subset \bigcup_{n\in \N}A_n$:
$$+\infty=\mathcal{H}^t(A_{n_0}) \leq \mathcal{H}^t(\bigcup_{n\in \N}A_n)$$
Therefore $\mathcal{H}^t(\bigcup_{n\in \N}A_n)=+\infty$.

The number s satisfies the property of Hausdorff dimension, hence:
$$\dim_{\mathcal{H}}(\bigcup_k A_k)=\sup \{\dim_{\mathcal{H}}(A_k):k\in \N\}$$
\end{proof}

\section{Application: Cantor Ternary Set}
\label{sec:applications}
We consider the interval $K_0=[0,1]$ with its middle third removed (excluding endpoints) $]\frac{1}{3},\frac{2}{3}[$.

\begin{figure}[H]
    \centering
    \includegraphics[scale=1]{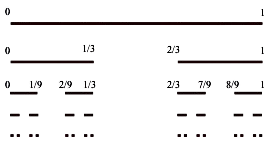}
    \caption{First five construction stages of the Cantor set}
    \label{fig:cantor-set}
\end{figure}

This results in the disjoint union of two closed intervals $K_1=[0,\frac{1}{3}]\cup[\frac{2}{3},1]$. Removing the middle third from each interval gives $K_2=[0,\frac{1}{9}]\cup[\frac{2}{9},\frac{1}{3}]\cup[\frac{2}{3},\frac{7}{9}]\cup[\frac{8}{9},1]$, a disjoint union of four closed intervals (Figure \ref{fig:cantor-set}).

Repeating this process n times yields $K_n$ - a union of $2^n$ disjoint closed intervals denoted $J_{n,k}$. We construct $K_{n+1}$ by removing middle thirds from each interval, producing a decreasing sequence $(K_n)_{n\in\N}$ of compact sets.

The intersection $K=\bigcap_{n\in\N} K_n$ is called the Cantor ternary set.

\textbf{Properties of K}:
\begin{itemize}
    \item \textbf{i)} K is a non-empty compact set.\\
As K is the intersection of a decreasing sequence of non-empty compact sets.
\item \textbf{ii)} K is uncountable.\\
Proof by contradiction: Assume there exists a bijection $f:\N^*\longrightarrow K$.\\
Let $I_1$ be the connected component of $K_1$ not containing $f(1)$, $I_2$ the connected component of $K_2$ such that $I_2\subset I_1$ and not containing $f(2)$, and generally, given $I_n$, let $I_{n+1}$ be the connected component of $K_{n+1}$ such that $I_{n+1}\subset I_n$ and not containing $f(n+1)$.\\
This produces a nested sequence $(I_n)_{n\in\N}$ of closed intervals satisfying:
$$ \bigcap_{n\geq 1} I_n\subset \bigcap_{n\geq 1}K_n=K $$
and
$$(\bigcap_{n\geq 1} I_n)\cap K=(\bigcap_{n\geq 1} I_n)\cap f(\N^*)=\varnothing $$
Thus $\bigcap_{n\geq 1} I_n=\varnothing$, contradicting the Cantor intersection theorem.

\item \textbf{iii)} K has Lebesgue measure zero.\\
As K is the intersection of nested compact sets $K_n$, by measure properties we have $\lambda(K)=\lim_{n\to\infty} \lambda(K_n)$, where $\lambda$ denotes Lebesgue measure.\\
Since each $K_n$ consists of $2^n$ disjoint intervals of length $\frac{1}{3^n}$, we have $\lambda(K_n)=\frac{2^n}{3^n}=(\frac{2}{3})^n$.\\
Taking the limit:
$$ \lambda(K)=\lim_{n\to\infty} (\frac{2}{3})^n=0$$
\end{itemize}

\begin{Th}
\label{th:4}
For the Cantor ternary set K and $s=\frac{\ln 2}{\ln 3}$:
$$\frac{1}{2} \leq\mathcal{H}^s(K)\leq1\;\;\text{and}\;\; \dim_{\mathcal{H}}(K)=s.$$
\end{Th}
\vspace{-2cm}
\begin{proof}[Proof of Theorem \ref{th:4}: ]
Let $K$ be the Cantor set and $s=\frac{\ln 2}{\ln 3}$.
\begin{itemize}
    \item \textbf{a)} Upper bound $\mathcal{H}^s(K)\leq 1$:\\
For any $n\in\N$, $K_n=\bigcup_{k=1}^{2^n}J_{n,k}$ with $diam(J_{n,i})=\frac{1}{3^n}$.\\
Since $\{J_{n,i}\}_{1\leq i\leq 2^n}\in\mathcal{R}_{3^{-n}}(K)$:
\begin{align*}
\mathcal{H}^s(K)&\leq \sum_{i=1}^{2^n}\left(\frac{1}{3^n}\right)^s\\
&= 2^n\left(3^{-n}\right)^{\frac{\ln 2}{\ln 3}}\\
&= 2^n2^{-n} = 1
\end{align*}

\item \textbf{b)} Lower bound $\mathcal{H}^s(K)\geq\frac{1}{2}$:\\
Preliminary: For any interval $I\subset\R$ with $0<diam(I)<\frac{1}{3}$, there exists unique $k\geq2$ such that:
$$3^{-k}\leq diam(I)\leq 3^{-(k-1)}$$
By construction, $I$ intersects at most one interval of $K_{k-1}$ and at most $2^{m-k+1}$ intervals of $K_m$ for $m>k$.

Main proof: Let $a\geq1$, $0<\delta<\frac{1}{3a}$, and ${U_i}\in\mathcal{R}_\delta(K)$. Take open intervals $T_i\supset U_i$ with $diam(T_i)=a\cdot diam(U_i)$. By compactness, extract finite subcover ${T_1,...,T_n}$. For $m$ satisfying:
$$3^{-(m-1)} < diam(T_i)\;\;\forall i=1,...,n$$
Each $T_i$ contains at most $2^{m+1}(diam(T_i))^s$ intervals of $K_m$. Since all $2^m$ intervals of $K_m$ must be covered:
$$2^m \leq \sum_{i=1}^n 2^{m+1}(a\cdot diam(U_i))^s$$
Thus:
$$\frac{1}{2a^s} \leq \sum_{i=1}^\infty diam(U_i))^s$$
Taking $a\to1$ and $\delta\to0$ gives $\mathcal{H}^s(K)\geq\frac{1}{2}$.

Conclusion: $\frac{1}{2}\leq\mathcal{H}^s(K)\leq1$ implies $\dim_{\mathcal{H}}(K)=s$.
\end{itemize}

\end{proof}
\bibliography{reference.bib}
\bibliographystyle{plainnat}

\end{document}